\documentstyle[amscd,amssymb,verbatim]{amsart}
\renewcommand{\ker}{\operatorname{ker}}
\newcommand{\im}{\operatorname{im}}
\newcommand{\D}{{\cal D}}
\numberwithin{equation}{subsection}
      
\newcommand{\lrar}[1]{\begin{picture}(50,10)(-25,-5)
\put(-25,0){\vector(1,0){50}}
\put(0,5){\makebox(0,0)[b]{\mbox{$#1$}}}
\end{picture}}

\newcommand{\ldar}[1]{\begin{picture}(10,50)(-5,-25)
\put(0,25){\vector(0,-1){50}}
\put(5,0){\mbox{$#1$}}
\end{picture}}

\newcommand{\luar}[1]{\begin{picture}(10,50)(-5,-25)
\put(0,-25){\vector(0,1){50}}
\put(5,0){\mbox{$#1$}}
\end{picture}}

\newcommand{\ldrar}[1]{\begin{picture}(50,50)(-25,-25)
\put(-25,25){\vector(1,-1){50}}
\put(5,0){\mbox{$#1$}}
\end{picture}}

\newcommand{\lular}[1]{\begin{picture}(50,50)(-25,-25)
\put(25,-25){\vector(-1,1){50}}
\put(5,0){\mbox{$#1$}}
\end{picture}}

\newcommand{\lurar}[1]{\begin{picture}(50,50)(-25,-25)
\put(-25,-25){\vector(1,1){50}}
\put(5,-5){\mbox{$#1$}}
\end{picture}}

\newtheorem{thm}{Theorem}[section]
\newtheorem{prop}[thm]{Proposition}
\newtheorem{lem}[thm]{Lemma}
\newtheorem{cor}[thm]{Corollary}
\newenvironment{rem}{\vspace{3mm}\noindent
{\bf Remark.}}{\vspace{3mm}}
\newcommand{\Pf}{\noindent {\it Proof}}
\newcommand{\ed}{\qed\vspace{3mm}}
\newcommand{\id}{\operatorname{id}}
\newcommand{\ov}{\overline}
\renewcommand{\Im}{\operatorname{Im}}
\newcommand{\ra}{\rightarrow}

\newcommand{\PP}{{\cal P}}
\newcommand{\LL}{{\cal L}}
\renewcommand{\O}{{\cal O}}
\newcommand{\Hom}{\operatorname{Hom}}
\newcommand{\Ext}{\operatorname{Ext}}
\newcommand{\Res}{\operatorname{Res}}

\renewcommand{\a}{\alpha}
\renewcommand{\b}{\beta}
\newcommand{\om}{\omega}
\newcommand{\De}{\Delta}
\newcommand{\la}{\lambda}
\newcommand{\th}{\theta}
\newcommand{\C}{{\Bbb C}}
\newcommand{\R}{{\Bbb R}}
\newcommand{\Z}{{\Bbb Z}}
\newcommand{\Q}{{\Bbb Q}}
\newcommand{\La}{\Lambda}
\newcommand{\Ga}{\Gamma}
\newcommand{\e}{\operatorname{e}}
\newcommand{\wt}{\widetilde}
\newcommand{\ot}{\otimes}
\newcommand{\sign}{\operatorname{sign}}
\newcommand{\sub}{\subset}
\newcommand{\eps}{\varepsilon}

\title{Massey and Fukaya products on elliptic curves}

\begin{document}
\maketitle

\centerline{A. Polishchuk\footnote{E-mail: apolish@@math.harvard.edu}}

\bigskip

\centerline{Department of Mathematics}
\centerline{Harvard University, Cambridge, MA 02138}       

\vspace{6mm}
     
This note is a complement to \cite{PZ}.
One of its goals is to show that some higher Massey products
on an elliptic curve can be computed as higher compositions in
Fukaya category of the dual symplectic torus in accordance with the
homological mirror conjecture of M.~Kontsevich \cite{Kon}.
Namely, we consider triple Massey products of very simple type which
are uniquely defined,
compute them in terms of theta-functions and compare the result with
the series one obtains in Fukaya picture. The identity we get in this
way was first discovered by Kronecker.
     
An interesting phenomenon is that although Massey products on an
elliptic curve
are partially defined and multivaled, one always has the
corresponding univalued Fukaya product. Thus, the equivalence with
Fukaya category equips the derived
category  of coherent sheaves on an elliptic curve with some
additional structure.
In order to understand this structure we study the relation
between the triple products $m_3$ and the triangulated  structure.
It  turns  out  that  using  $m_3$  one  can define homotopy
operators on cohomological  long  exact  sequences
associated with a generic distinguished triangle. Furthermore,
higher products $m_k$ with $k\ge 4$ define higher homotopy operators on these
exact sequences. It seems that most of the $A_{\infty}$-structure 
can be recovered from these homotopy operators.

The rest of the paper is devoted to explicit computations of
higher compositions $m_3$ in Fukaya  category of a torus
corresponding to four lines with rational slopes.
It turns out that the answer is given in terms of theta series associated with
not necessarily definite quadratic forms on rank-2 lattices.
Such series were introduced by L. G\"ottsche and D. Zagier in
\cite{GZ}. The idea is that when the quadratic form on the lattice is
indefinite one has to restrict the summation over the lattice to the cone
where the form is positive (introducing signs for different connected
components of the cone). In the case when one takes the maximal cone
on which the quadratic form is positive such a series is especially nice:
as shown in \cite{GZ} it is a Jacoby form.
In the case when the quadratic form factors over $\Q$
into product of two linear forms, indefinite theta series for
arbitrary cones can be expressed via the function
$$\kappa(y,x;\tau)=\sum_{n\in\Z}\frac{\exp(2\pi i (\tau n^2/2+nx))}
{\exp(2\pi i n\tau)-\exp(2\pi i y)}$$
where $\tau$ is in the upper half-plane. The latter function was
introduced (with slightly different notation) by M.~P.~Appell
in his study of doubly-periodic functions of the third kind in \cite{Ap}.
We write explicitly the $A_{\infty}$-constraints
between $m_2$ and $m_3$ in Fukaya category of a torus as
an identity between indefinite theta series. In particular, we
recover some
non-trivial identities involving $\kappa$ and theta-functions.
 
Our computation of univalued triple Massey
products on elliptic curve can be generalized to the case of
higher  genus  curves.  The answer is always given as certain
ratio of theta-functions.
We expect that these products can be
compared with Fukaya compositions on the symplectic
torus, which is mirror dual to the Jacobian of a curve.
This generalization will be considered in a future paper.

Throughout this paper we use the notation
$\e(z)=\exp(2\pi i z)$.
         
\vspace{3mm}

\noindent
{\it Acknowledgment}. I am grateful to  M.~Kontsevich  for
explaning to me the signs in Fukaya composition and to E.~Zaslow
for correcting some errors. This work was partially supported
by NSF grant.

\section{Triple Massey products}
       
\subsection{General construction}\label{gen-con}(cf. \cite{GM})
Let $X_1\stackrel{f_1}{\ra} X_2\stackrel{f_2}{\ra} X_3\stackrel{f_3}{\ra}X_4$
be a sequence of morphisms in a triangulated category such that
$f_2\circ f_1=0$ and $f_3\circ f_2=0$. Then one can construct the
subset of elements $MP(f_1,f_2,f_3)$
in $\Hom^{-1}(X_1,X_4)$ which is a coset by the sum of the images
of the maps
$$\Hom^{-1}(X_2,X_4)\ra\Hom^{-1}(X_1,X_4):g\mapsto g\circ f_1,$$
$$\Hom^{-1}(X_1,X_3)\ra\Hom^{-1}(X_1,X_4):h\mapsto f_3\circ h.$$
Namely, $k\in MP(f_1,f_2,f_3)$ iff $k=v\circ u$ and
the following diagram is commutative
\begin{equation}
\begin{array}{cccc}
X_2&\lrar{f_2}&X_3\\
\luar{f_1}&\lular{[1]}&\luar{}&\ldrar{f_3}\\
X_1&\lrar{u}&Y&\lrar{v}X_4
\end{array}
\end{equation}
where $X_2\ra X_3\ra Y\ra X_2[1]$ is an  exact triangle.
In particular, if $\Hom^i(X_1,X_3)=\Hom^i(X_2,X_4)=0$ for $i=-1,0$
then the compositions $f_2\circ f_1$ and $f_3\circ f_2$ are always zero,
hence, the Massey
product $MP(f_1,f_2,f_3)$
contains one element, and we obtain the linear map
$$m_3:\Hom(X_1,X_2)\otimes\Hom(X_2,X_3)\otimes\Hom(X_3,X_4)\ra
\Hom^{-1}(X_1,X_4).$$
Here is a more concrete description of this map.
Consider the exact triangle
$$K\ra\Hom(X_2,X_3)\otimes X_2\ra X_3\ra K[1]$$
Then  our  assumptions  imply  that  the  following  natural
maps are isomorphisms:
$$\alpha:\Hom(X_1,K)\ra\Hom(X_1,X_2)\otimes\Hom(X_2,X_3),$$
$$\beta:\Hom^{-1}(K,X_4)\ra\Hom(X_3,X_4).$$
Now $m_3$ is equal to the following composition
\begin{align*}
&\Hom(X_1,X_2)\otimes\Hom(X_2,X_3)\otimes\Hom(X_3,X_4)
\lrar{\a^{-1}\otimes\b^{-1}}\\
&\Hom(X_1,K)\otimes\Hom^{-1}(K,X_4)\ra\Hom^{-1}(X_1,X_4)
\end{align*}
where the last arrow is the natural
composition map.
     
\subsection{Some Massey products on elliptic curve}\label{simple}
       
Let us consider the simplest example
of triple Massey product on elliptic curve $E$ over a field $k$.
Namely, we want to descibe the Massey product
$$\Hom(\O,\O_{x_0})\otimes\Hom(\O_{x_0},\LL[1])\otimes
\Hom(\LL[1],\O_x[1])\ra\Hom(\O,\O_x)$$
where $x\neq x_0$, $\LL\neq\O$, $\deg\LL=0$. If $\LL|_{x_0}$
is trivialized then using the Serre duality the source of
this arrow can be identified with $(\LL|_x\otimes\om)^*$
where $\om$ is the stalk of the canonical bundle of $E$ at zero,
while the target is $k$. The dual map gives a canonical element
$s(\LL,x)\in\LL|_x\otimes\om$. Note that $\LL|_x$ is a stalk of the Poincar\'e
line bundle $\PP$ on $E\times\hat{E}$ where $\hat{E}$ is the dual elliptic
curve. It is easy to see the above Massey product is the value at
the point $(\LL,x)$ of the canonical rational section of
$\PP\otimes\om$ (with poles at $x=x_0$ and $\LL=\O$).
            
Now assume  that  $k=\C$  and  the  elliptic  curve  $E$  is
$\C/\Ga$ where $\Ga=\Ga_{\tau}=\Z+\Z\tau$,
$\tau$ is in the upper-half plane. We want to express the above Massey product
in terms of theta-function. Let us denote by $\pi:\C\ra E$ the canonical
projection.
Consider the line bundle $L$ on $E$ (equipped with a trivialization
of its pull-back to $\C/\Z$) such that the classical
theta-function
$$\theta(z,\tau)=\sum_n\e(\tau n^2/2+ nz)$$
is a section of $L$ (as a function of $z$),
so $L\simeq\O(\xi)$ where $\xi=(\tau+1)/2$.
Now set $\LL=t^*_yL\otimes L^{-1}$ where $y\in\C$ is not a lattice
point. We set $x_0=\pi(0)$ and fix a lifting of the second point $x$ to
$\C$ (abusing the notation we denote this lifting by $x\in\C$).
Then the trivialization of $\pi^*L$ induces trivializations
of $\LL_{x_0}$ and $\LL_x$, so we have canonical generators
$f_1\in\Hom(\O,\O_{x_0})$, $f_2\in\Hom(\LL,\O_{x_0})^*$ and
$f_3\in\Hom(\LL,\O_x)$. For a non-zero global holomorphic
1-form $\a\in H^0(E,\om)$ we can consider the corresponding isomorphism
of functors
\begin{equation}\label{duality}
S_{\a}:\Hom(A,B)^*\ra\Hom^1(B,A)
\end{equation}
derived from the Serre duality. Then we have an element
$S_{\a}(f_2)\in\Hom(\O_{x_0},\LL[1])$, and we can consider the triple
Massey product $MP(f_1,S_{\a}(f_2),f_3)$.

\begin{lem}\label{1form}
Let $f\in\Hom(\O,\O_{\xi})$ be the canonical generator,
let $\a=\th'(\frac{\tau+1}{2})dz$. Then the element
$S_{\a}(f)\in\Ext^1(\O_{\xi},\O)$ is represented by the
extension
\begin{equation}
\label{ext1}
0\ra \O\stackrel{\th}{\ra}L\ra\O_{\xi}\ra 0.
\end{equation}
\end{lem}

\Pf . We have the canonical extension
\begin{equation}\label{ext2}
0\ra\om\ra\om(\xi)\stackrel{\Res}{\ra}\O_{\xi}\ra 0.
\end{equation}
Via  the  isomorphism  $\O\simeq\om$  induced  by  $\a$ this
extension represents $S_{\a}(f)$.  Now  we  claim  that  the
map
$$L\ra\om(\xi):s\mapsto\frac{s}{\th}\cdot\a$$
extends   to   the   isomorphism  between  (\ref{ext1})  and
(\ref{ext2}). Indeed, this follows from the fact that
$\Res_{\xi}(\frac{\a}{\th})=1$.
\ed

\begin{prop} Let $\a=\th'(\frac{\tau+1}{2})dz$. Then
\begin{equation}\label{mas_main}
MP(f_1,S_{\a}(f_2),f_3)=
\frac{\th(x+y+\xi)}{\th(x+\xi)\th(y+\xi)}\cdot f_x
\end{equation}
where $f_x\in\Hom(\O,\O_x)$ is the canonical generator,
$\xi=\frac{\tau+1}{2}$.
\end{prop}

\Pf . It follows from the above Lemma that the element
$S_{\a}(f_2)$ correspond to the extension
$$0\ra \LL\stackrel{t_{\xi}^*\th}{\ra}t_{y+\xi}^*L\ra\O_{x_0}\ra 0.$$
The recipe for computation of Massey products is:
first lift $f_1$ to the section $s$ of $t_{y+\xi}^*L$;
then find the morphism $g:t_{y+\xi}^*L\ra\O_x$ such that
the diagram
\begin{equation}
\begin{array}{ccc}
\LL &\lrar{t_{\xi}^*\th}& t_{y+\xi}^*L\\
&\ldrar{f_3}&\ldar{g}\\
&&\O_x
\end{array}
\end{equation}
commutes; and finally apply $g$ to $s$.
We have $s=\frac{t_{y+\xi}^*\th}{\th(y+\xi)}$,
$g$ is $\frac{1}{\th(x+\xi)}$ times the canonical generator of
$\Hom(t_{y+\xi}^*L,\O_x)$  (induced  by  the   trivialization   of
$\pi^*L$), so we get (\ref{mas_main}).
\ed

\section{Comparison with Fukaya composition}

\subsection{Triple Fukaya composition and triangulated structure}
In  \cite{PZ} we have constructed an equivalence between the
derived category of coherent sheaves on an elliptic curve $E$ and the
Fukaya  category  of  the  corresponding 2-dimensional torus
with complexified symplectic form.
An object of the latter category consists
of the following data: a geodesic circle $\La\in T=\R^2/\Z^2$,
an angle $\phi\in\R$ such that $\R\e(\phi)$ is parallel to
$\La$,  and  a  local  system $\LL$ on $\La$. More precisely, one can
consider formal direct sums of such objects.
Note that the change of $\phi$
by  $\phi+1$  corresponds  to  the  translation functor on the
derived category.
Morphisms from $(\La,\phi,\LL)$  to  $(\La',\phi',\LL')$  (where
$\La\neq  \La'$) can be non-zero only if $\phi<\phi'<\phi+1$. In
this case the $\Hom$-space is
$$\sum_{x\in \La\cap \La'}\LL_x^*\otimes\LL_x.$$
The composition  is  defined
using holomorphic triangles bounding three given geodesics
circles (see \cite{PZ}). There are also higher compositions
$m_k$, $k\ge 3$ which are defined using holomorphic $(k+1)$-gons
(see \cite{F-rec}). They satisfy $A_{\infty}$-axioms 
(with $m_1=0$) of which the following is an example:
\begin{eqnarray}\label{a-infty}
\nonumber m_3(m_2(a_1,a_2),a_3,a_4)\pm m_3(a_1,m_2(a_2,a_3),a_4)\pm
m_3(a_1,a_2,m_2(a_3,a_4))\\
\pm m_2(m_3(a_1,a_2,a_3),a_4)\pm m_2(a_1,m_3(a_2,a_3,a_4))=0
\end{eqnarray}
where   $a_1,\ldots,a_4$   is  the  sequence  of  composable
morphisms, the signs depend on degrees of $a_i$'s.
Note that the proof of $A_{\infty}$-identities for the Fukaya
compositions in our case is easy. Basically they follow from
additivity of the area of plane figures. For example, the usual
associativity of $m_2$ (which is equivalent to the addition formula for
theta function) can be proved by considering two ways of
cutting a non-convex quadrangle into two triangles and converting
this into the identity for the corresponding generating series.

The functor from the Fukaya category to the derived category
of coherent sheaves on elliptic curve is set up in such a way
that one gets objects of abelian category (i.e. complexes
concentrated in degree zero) for $\phi\in (-1/2,1/2]$.
The  functor is constructed on  objects  from  this
subcategory using theta-functions.
To obtain the entire derived category $\D^b(E)$ one has
to deal with $\Ext^1$.
We fix a generator $\a\in H^0(E,\om)$ and
use the corresponding duality isomorphism (\ref{duality})
(on Fukaya side the corresponding isomorphism is obvious).
Then we get a functor $\Phi_{\a}$ from the Fukaya category
to $\D^b(E)$.

In \cite{F-rec} Fukaya constructs the functor on $\Ext^1$ in
a different way. Namely, he starts  with  the  same  functor
$\Phi$
on abelian category and then proceeds as follows.
Let $a:X\ra Y$ be a morphism of degree 1 in Fukaya category.
To construct the corresponding morphism  $\Phi(a):\Phi(X)\ra
\Phi(Y)[1]$
one has to choose a resolution
$$0\ra           \Phi(Y)\stackrel{\Phi(d_0)}{\ra}
\Phi(Z_0)\stackrel{\Phi(d_1)}{\ra} \Phi(Z_1)
\stackrel{\Phi(d_2)}{\ra}\Phi(Z_2)$$
such that $\Hom^1(X, Z_0)=0$. Now consider
the triple composition $m_3(a,d_0,d_1)\in\Hom^0(X,Z_1)$.
Then          $A_{\infty}$-axioms         imply         that
$m_2(m_3(a,d_0,d_1),d_2)=0$, hence, $-\Phi(m_3(a,d_0,d_1))$
factors through $\ker \Phi(d_2)=\im \Phi(d_1)$, so it defines an
element in $\Hom^1(\Phi(X),\Phi(Y))$ which we take to be $\Phi(a)$
(we have changed the sign compare to Fukaya's definition in order
for Proposition \ref{ex-tr} below to be true).
In  fact,  in   the   above   construction   the   condition
$\Hom^1(X,Z_1)=0$ can be relaxed to the requirement that
$m_2(a,d_0)=0$. The independence on a choice of a resolution
is proven by providing an alternative definition via a resolution for
$\Phi(X)$.  Let  us denote by $\Phi$ the obtained functor from the
Fukaya category to $\D^b(E)$. Fukaya showed in  \cite{F-rec}
that $\Phi$ is an equivalence so it should coincide
with $\Phi_{\a}$ for certain 1-form $\a$.
To find $\a$ we  will  use  the  following characterizing property of the
functor $F$.

\begin{prop}\label{ex-tr} Let    $a\in\Hom(X,Y)$,     $b\in\Hom(Y,Z)$,
$c\in\Hom(Z,X[1])$ be morphisms in Fukaya category such that
$\Phi(a)$, $\Phi(b)$ and $\Phi(c)$ form an exact triangle. Assume
in addition that $\Hom(Y,X)=0$. Then $m_3(a,b,c)=\id_X$.
\end{prop}

\Pf . Let us use the above definition to compute $\Phi(c)$.
Namely, let us choose a resolution for $\Phi(Z)$ of the form
$$0\ra\Phi(Z)\stackrel{\Phi(d_1)}{\ra}\Phi(T_1)
\stackrel{\Phi(d_2)}{\ra}\Phi(T_2).$$
Then we have the induced resolution for $\Phi(X)$ which
fits into the commutative diagram
\begin{equation}
\begin{array}{cccc}
0\ra\Phi(X)\lrar{\Phi(a)} &\Phi(Y)&
\lrar{}&\Phi(T_1)
\lrar{\Phi(d_2)}\Phi(T_2)\\
&\ldar{\Phi(b)}&\lurar{\Phi(d_1)}\\
&\Phi(Z)
\end{array}
\end{equation}
Then by definition $\Phi(c)$ is induced by the morphism
$-\Phi(m_3(c,a,m_2(b,d_1)))$. Now the $A_{\infty}$-constraint
implies that
$$m_3(c,a,m_2(b,d_1))=m_2(m_3(c,a,b),d_1).$$
It follows that
$\Phi(c)$ is the composition of $-\Phi(m_3(c,a,b))$
with $\Phi(c)$, i.e. $m_2(m_3(c,a,b),c)=-c$.
Now again by $A_{\infty}$-axiom we deduce that
$$m_2(c,m_3(a,b,c))=-m_2(m_3(c,a,b),c)=c.$$
Since $\Hom(Y,X)=0$ this implies that $m_3(a,b,c)=\id_X$.
\ed

The following corollary can be found in Fukaya's paper \cite{F-rec}.
\begin{cor}\label{comp-Mas} Let $a\in\Hom(X,Y)$, $b\in\Hom(Y,Z)$,
$c\in\Hom(Z,X)$ be a triple of morphisms in Fukaya
category such that $b\circ a=0$ and $c\circ b=0$.
Then $\pm\Phi(m_3(a,b,c))\in MP(\Phi(a),\Phi(b),\Phi(c))$.
\end{cor}

\begin{cor}\label{homotopy} Let
$$X\stackrel{a}{\ra} Y\stackrel{b}{\ra} Z\stackrel{c}{\ra} X[1]$$
be a distinguished triangle in $\D^b(E)$ identified with
the Fukaya category using the functor $\Phi$. Assume that
$$\Hom^i(X,Y)=\Hom^i(Y,Z)=\Hom^{i+1}(Z,X)=0$$
for $i\neq 0$.
Then for every object $T\in\D^b(E)$ there is a canonical homotopy operator
$H_1$ on the corresponding long exact sequence of morphisms from $T$ to the
above triangle, i.e. $H_1\partial+\partial H_1=\id$ where $\partial$ is the
differential in the latter exact sequence. Namely, for
$f\in\Hom(T,X)$ one should take $H_1(f)=\pm m_3(f,a,b)\in\Hom(T,Z[-1])$, etc.
Furthermore, one can define higher homotopy operators by setting
$H_2(f)=\pm m_4(f,a,b,c)$, $H_3(f)=\pm m_5(f,a,b,c,a)$, etc.,
such that $H_1^2=H_2\partial+\partial H_2$, 
$H_1H_2-H_2H_1=H_3\partial+\partial H_3$, etc.
Similarly, there are canonical homotopy operators $H'_i$
for  the  long exact sequence of morphisms from the triangle
to $T$. 
\end{cor}

The proofs of both corollaries are easy exercises
in applying $A_{\infty}$-axioms which
we leave to the reader. In particular, for the second corollary one has
to use that fact that higher products containing an identity morphism
vanish.

\begin{prop}\label{pi} One has $\Phi=\Phi_{2\pi i dz}$.
\end{prop}

\Pf  .  Let  us  consider the following three of objects in
Fukaya category: $\La_1=(t,0)$ with trivial local system,
$\La_2=(t,t)$ with trivial local system,
$\La_3=(1/2,t)$ with the connection $\pi i dt$.
Then we have  canonical  morphisms $e_i$  from  $\La_i$  to
$\La_{i+1}$    such    that    $\deg(e_1)=\deg(e_2)=0$   and
$\deg(e_3)=1$. It follows easily from Lemma \ref{1form}
that the morphisms
$\Phi(e_1)$,  $\Phi(e_2)$  and  $\Phi_{\th'(\xi)dz}(e_3)$  form  an  exact
triangle
$$\O\ra L\ra \O_{\xi}\ra\O[1]$$
where $\xi=\frac{\tau+1}{2}$.
On the other hand, it is easy to compute that
$$m_3(e_1,e_2,e_3)=\frac{\th'(\xi)}{2\pi i}.$$
It follows that
$m_3(e_1,e_2,\frac{2\pi i}{\th'(\xi)}e_3)=1$ while
$\Phi(e_1)$,  $\Phi(e_2)$  and
$\Phi_{2\pi i dz}(\frac{2\pi i}{\th'(\xi)}e_3)$ form
an exact triangle, hence $\Phi_{2\pi i dz}=\Phi$.
\ed

\subsection{Fukaya series}
Let   us   consider   the  triple  Fukaya  composition  that
corresponds to the Massey product considered in
\ref{simple}.
The corresponding four objects of the Fukaya category of a torus
are:
$\La_1=(t,0)$ with trivial connection,
$\La_2=(0,t)$ with trivial connection,
$\La_3=(t,-\a_1)$ with the connection $(-2\pi i\b_1)dx$,
$\La_4=(-\a_2,t)$ with the connection $2\pi i\b_2dy$.
Here $\a_1,\a_2,\b_1,\b_2$ are real numbers.
There is an essentially unique choice of logarithms of slopes for $\La_i$
such that $\Hom^0(\La_i,\La_{i+1})\neq 0$ and with such a choice one
has $\Hom^{-1}(\La_1,\La_4)\neq 0$. Moreover, all these spaces are
one-dimensional so the corresponding Fukaya composition $m_3$ is just
a number. This composition is defined only if $\a_1$ and $\a_2$ are
not integers. Then it is given by the following series
$$\sum_{(\a_1+m)(\a_2+n)>0}\sign(\a_1+m)
\e(\tau(\a_1+m)(\a_2+n)+(\a_1+m)\b_2+(\a_2+n)\b_1)$$
where we denote $\sign(t)=1$ for $t>0$, $\sign(t)=-1$ for $t<0$,
the sum is over integers $m$ and $n$ subject to the
condition that $(\a_1+m)$ has the same sign as $(\a_2+n)$,
The restriction on $m$ and $n$ is imposed by the condition
for maps from the disk to be holomorphic in the definition of
Fukaya composition while the sign comes from the canonical orientaion of
the corresponding moduli space.
We can write the above expression in the form
$$\e(\tau\a_1\a_2+\a_1\b_2+\a_2\b_1)f(\a_1\tau+\b_1,
\a_2\tau+\b_2;\tau)$$
where
$$f(z_1,z_2;\tau)=\sum_{(\a(z_1)+m)(\a(z_2)+n)>0}
\sign(\a(z_1)+m)
\e(\tau mn+nz_1+mz_2)$$
is a holomorphic function of $z_1$ and $z_2$ defined for
$\Im z_i\not\in\Z(\Im\tau)$, where $\a(z)=\Im(z)/\Im(\tau)$.

\subsection{The function $f(z_1,z_2,\tau)$}
The function $f(z_1,z_2;\tau)$ is well-known (cf. \cite{Kr},
\cite{Weil},\cite{Z}).
It extends to a meromorphic function in $z_1, z_2$ with
poles at the lattice points $z_1\in\Ga_{\tau}$ or
$z_2\in\Ga_{\tau}$ where $\Ga_{\tau}=\Z+\Z\tau$, and satisfies
the following identities:
\begin{eqnarray}\label{t1}
f(z_2,z_1;\tau)=f(z_1,z_2;\tau),\\
\label{t2}
f(z_1+m+n\tau,z_2;\tau)=\e(- n z_2)f(z_1,z_2;\tau),\\
\label{t3}
f(z_1,z_2+m+n\tau;\tau)=\e(- n z_1)f(z_1,z_2;\tau).
\end{eqnarray}

Using this quasi-periodicity properties of $f$ it is easy to
derive the following
identity (\ref{id}) which was first discovered by Kronecker \cite{Kr}
(see also \cite{Weil}, ch. VIII):
\begin{equation}\label{id}
f(z_1,z_2;\tau)=\frac{\th'((\tau+1)/2,\tau)}{2\pi i}
\cdot
\frac{\theta(z_1+z_2-(\tau+1)/2,\tau)}
{\theta(z_1-(\tau+1)/2,\tau)\theta(z_2-(\tau+1)/2,\tau)}
\end{equation}
Note that
$\frac{\th'((\tau+1)/2,\tau)}{2\pi  i}=\prod_{n\ge 1}(1-q^n)^3$
where $q=\e(\tau)$.

It is clear from the definition that
$$f(z_1,z_2;\tau+1)=f(z_1,z_2;\tau).$$
Now using the identity (\ref{id})
one can easily deduce from the functional equation for theta function that
$f(z_1,z_2;\tau)$ satisfies the functional equation of the form
$$f(z_1/\tau,z_2/\tau;-\tau^{-1})=\zeta\cdot\tau\cdot
\e(z_1z_2/\tau)f(z_1,z_2;\tau)$$
where $\zeta$ is a root of unity. Using the property
$$f(-z_1,-z_2;\tau)=-f(z_1,z_2;\tau)$$
we immediately conclude that $\zeta^2=1$, so $\zeta=\pm 1$.
Finally substituting $z_1=(\tau+1)/2$, $z_2=\tau/2$ and
looking at the sign of both sides when $\tau=it$, $t\in\R$ and
$t\ra+\infty$ we find that $\zeta=1$. So the functional equation for
$f$ becomes
$$f(z_1/\tau,z_2/\tau;-\tau^{-1})=\tau\cdot
\e(z_1z_2/\tau)f(z_1,z_2;\tau).$$
In fact, $f(z_1,z_2,\tau)$ is a meromorphic Jacobi form of weight
1 for the lattice $\Z^2$ with the quadratic form $Q(m,n)=mn$
(cf. \cite{GZ}).

This equation can also be derived from the representation of $f$
in the following form (cf. \cite{Weil}, ch.VIII):
$$f(z_1,z_2;\tau)=-\frac{\e(-\a(z_2)z_1)}{2\pi i}
\sideset{}{_e}{\sum} \frac{\chi(w)}{z_1+w}$$
provided that $0<\a(z_i)<1$ for $i=1,2$. Here $\sum_e$ denotes the
Eisenstein summation over the lattice $\Ga_{\tau}$,
$\chi$ is the character of $\Ga_{\tau}$ such that
$\chi(1)=\e(-\a(z_2))$, $\chi(\tau)=\e(z_2-\a(z_2)\tau)$.
 
\subsection{Comparison}
Now  we  want  to  interpret  the  identity (\ref{id}) as an
equality of the triple Fukaya product with the corresponding
Massey product asserted in Corollary \ref{comp-Mas}.
Using the explicit construction of the functor $\Phi$ (cf. \cite{PZ})
we compute that
$\Phi(\La_1)=\O$, $\Phi(\La_2)=\O_{x_0}$,
$\Phi(\La_3)=t_y^*L\otimes L^{-1}$ where $y=\a_1\tau+\b_1$,
$\Phi(\La_4)=\O_x$ where $x=\a_2\tau+\b_2$.
Let  $e_i\in\Hom^*(\La_i,\La_{i+1})$,  $i=1,2,3$, be canonical
generators, where $\deg(e_1)=\deg(e_3)=0$, $\deg(e_2)=1$.
Then using the notation of \ref{simple} we have
$\Phi(e_1)=f_1$, $\Phi_{\a}(e_2)=S_{\a}(f_2)$, and
$\Phi(e_3)=\e(\tau\a_1\a_2+\a_1\b_2+\a_2\b_1)f_3.$
Also for a canonical generator $e\in\Hom^0(\La_1,\La_4)$  we
have $\Phi(e)=f_x$.
Now for $\a=\th'(\xi)$, where $\xi=(\tau+1)/2$ we derive from
(\ref{mas_main}) that
$$MP(\Phi(e_1),\Phi_{\th'(\xi)dz}(e_2),\Phi(e_3))=
\e(\tau\a_1\a_2+\a_1\b_2+\a_2\b_1)\cdot
\frac{\th(x+y+\xi)}{\th(x+\xi)\th(y+\xi)}\cdot f_x$$
Thus using the above computation of $m_3(e_1,e_2,e_3)$
and Proposition \ref{pi} we get
\begin{align*}
&\e(\tau\a_1\a_2+\a_1\b_2+\a_2\b_1)f(x,y)\cdot f_x=
\Phi(m_3(e_1,e_2,e_3))=\\
&MP(\Phi(e_1),\Phi_{2\pi i dz}(e_2),\Phi(e_3))=
\frac{\th'(\xi)}{2\pi i}
MP(\Phi(e_1),\Phi_{\th'(\xi)dz}(e_2),\Phi(e_3)).
\end{align*}
Using the above expression for the Massey product in the RHS
we obtain the identity (\ref{id}).
                             
\begin{rem} It is not hard to see  that  Fukaya  triple
products involving four lines forming any parallelogram with
sides of rational slopes are expressed via the function $f$
and the equality with the corresponding (univalued) Massey
products follows from the identity (\ref{id}).
\end{rem}

\section{Higher compositions in Fukaya category}
   
\subsection{Trapezoid compositions}
Now we are going to consider some compositions
$m_3$  in  Fukaya  category  of  a  torus  such   that   the
corresponding triple Massey products on elliptic curve
are not well-defined. Namely, consider four lagrangians:
$\La_1=(t,-t)$ and $\La_4=(t,0)$ with trivial connections,
$\La_2=(t,-\a_2)$ with the connection $-2\pi i\b_2 dx$, and
$\La_3=(-\a_1,t)$ with the connection $2\pi i\b_1 dy$,
where $\a_i$, $\b_i$ are real numbers, $\a_i\not\in\Z$.
There is a natural choice of logarithms of slopes so that
$\Hom^0(\La_i,\La_{i+1})\simeq\C$ and
$\Hom^{-1}(\La_1,\La_4)\simeq\C$.
Note that $\Hom^0(\La_1,\La_3)\neq 0$, so the corresponding
triple Massey product on elliptic curve is not defined.
The Fukaya composition
$$m_3:\Hom(\La_1,\La_2)\ot\Hom(\La_2,\La_3)\ot
\Hom(\La_3,\La_4)\ra\Hom^{-1}(\La_1,\La_4)$$
is just the number given by the series
$$\sum_{(n+\a_1)(m+\a_2)>0}
\sign(m+\a_2)\e((n+\a_1+\frac{m+\a_2}{2})(m+\a_2)\tau+
(m+\a_2)\b_1+(m+\a_2+n+\a_1)\b_2).$$
   
Let us define
$$g(z_1,z_2;\tau)=\sum_{(n+\a(z_1))(m+\a(z_2))>0}
\sign(m+\a(z_2))\e((n+m/2)m\tau+mz_1+(m+n)z_2)$$
where as before $\a(z)=\Im(z)/\Im(\tau)$,
$\a(z_1)\not\in\Z$, $\a(z_2)\not\in\Z$.
Then the above composition is equal to
$$\e((\a_1+\a_2/2)\a_2\tau+\a_2\b_1+(\a_1+\a_2)\b_2)
g(\a_1\tau+\b_1,\a_2\tau+\b_2;\tau).$$
   
\subsection{Properties of $g(z_1,z_2;\tau)$}
The function $g$ is holomorphic for $\a(z_i)\not\in\Z$,
$i=1,2$, and
satisfies the following quasi-periodicity identities
\begin{eqnarray}
\label{t4}
g(z_1+m+n\tau,z_2;\tau)=\e(-n z_2)g(z_1,z_2;\tau),\\
\label{t5}
g(z_1,z_2+m+n\tau;\tau)=\e(-n^2\tau/2-n(z_1+z_2))g(z_1,z_2;\tau).
\end{eqnarray}
In other words, $g$ can be considered as a holomorphic
section of a line bundle on $E^2$ over the open subset $(E\setminus S)^2$
where $S=\R/\Z\subset\C/\Ga_{\tau}$.
However, it doesn't extend to a meromorphic section on
$E^2$. On the other hand, let us denote by $g_0$ the restriction of
$g$ to the region $0<\a(z_1)<1$, $\a(z_2)\not\in\Z$. Then
we claim that $g_0$ extends to a meromorphic function on $\C^2$
with poles of order 1 at $z_2\in\Ga_{\tau}$. Indeed,
if we sum first over $n$ in the series defining $g$ we get
$$g_0(z_1,z_2;\tau)=
\sum_{m\in\Z}\frac{\e(m^2\tau/2+m(z_1+z_2))}{1-\e(m\tau+z_2)}.$$
The latter series clearly satisfies the properties we claimed.
However, $g_0$ lacks the quasi-periodicity of
$g$. More precisely, it is easy to see that $g$ and $g_0$ are related
as follows:
$$g_0(z_1,z_2;\tau)-g(z_1,z_2;\tau)=p(z_1,\tau)\th(z_1+z_2;\tau)$$
where $p(z,\tau)$ is the following piecewise polynomial function of $\e(z_1)$:
$$p(z,\tau)= \cases -\sum_{0<n\le \a(z)}\e(-n^2\tau/2+nz_1), &\a(z)\ge 0,\\
\sum_{\a(z)<n\le 0} \e(-n^2\tau/2+nz_1), &\a(z)<0 \endcases$$

Let us consider the following series
\begin{equation}
\kappa(y,x;\tau)=\sum_{n\in\Z}\frac{\e(n^2/2\tau+nx)}
{\e(n\tau)-\e(y)}.
\end{equation}
This function is holomorphic for $y\not\in\Ga_{\tau}$
and satisfies the difference equation
$$\kappa(y,x+m+\tau;\tau)=\e(y)\kappa(y,x;\tau)+
\theta(x,\tau)$$
where $m\in\Z$.
Then we have
\begin{equation}\label{g0}
g_0(z_1,z_2;\tau)=-\e(-z_2)\kappa(-z_2,z_1+z_2;\tau)=
\kappa(z_2,\tau-z_1-z_2;\tau).
\end{equation}

\begin{rem} The function $\kappa$ and its derivatives were used by M.~P.~Appell
to represent an arbitrary doubly-periodic function of the third kind as
a sum of simple elements (cf. \cite{Ap},\cite{H}).
On the other hand, $\kappa$ can be expressed via the bilateral basic
hypergeometric series. Namely, using notation of \cite{GR}
we have
$$\kappa(y,x;\tau)=(1-\e(-y))^{-1}\cdot
\sideset{_1}{_2}{\psi}(\e(-y);0,\e(\tau-y);\e(\tau),
\e(x+(\tau+1)/2)).$$
\end{rem}

\subsection{Trapezoid Massey products}
The geometric meaning of the function $\kappa$ above is the following:
the pair $(\kappa(y,x;\tau),\theta(x,\tau))$ determines a global section
of a rank 2 bundle $F_y$ on $E$ which is a non-trivial extension of $L$
by the line bundle of degree 0 corresponding to $y$.
Namely, the bundle $F_y$ is defined as follows
$$F_y=\C^*\times\C^2/(z,v)\mapsto (z\cdot\e(\tau),A_y(z)v)$$
where
$$A_y(z)=\left(  \matrix \e(y) & 1 \\ 0 & \e(-\tau/2)z^{-1} \endmatrix
\right).$$
We claim that some of triple Massey products corresponding to the
trapezoid Fukaya products considered above are univalued and are
also expressed via $\kappa$. Namely, consider the
triple Massey product for $X_1=\O$, $X_2=L$, $X_3=\LL[1]$, $X_4=\O_{\xi}[1]$
where $\LL$ is a non-trivial line bundle of degree 0, $\xi=(\tau+1)/2$. Then
$\Hom^*(X_1,X_3)=0$ and $\Hom^0(X_2,X_4)=0$. Furthermore, the composition
map
$$\Hom(X_1,X_2)\otimes\Hom^{-1}(X_2,X_4)\ra\Hom^{-1}(X_1,X_4)$$
is zero since the unique section of $L$ vanishes at $\xi$. It follows
that the Massey product is well-defined and univalued in this case.
To compute it one should include the non-zero morphism $X_2\ra X_3$
into an exact triangle
$$X_2\ra X_3\ra C\ra X_2[1]$$
then lift the morphisms $X_1\ra X_2$ and $X_3\ra X_4$ to morphisms
$X_1\ra C[-1]$ and $C\ra X_4$ and compose the obtained two morphisms.
In our case one starts with a global section $s:\O\ra L$, then
$s$ can be lifted canonically to a section $\wt{s}$ of $F_y$ where
$\LL\simeq t_yL^{-1}\otimes L$. Now one chooses a splitting
$r:(F_y)_{\xi}\ra\C$ of the embedding $\C\simeq \LL_{\xi}\ra (F_y)$
and applies $r$ to $\wt{s}(\xi)$. The result doesn't depend
on a choice of splitting at $\xi$ since $s(\xi)=0$. More concretely,
for $s=\theta(x,\tau)$ we have $\wt{s}=\kappa(y,x;\tau)$, hence, the
above Massey product is given by $\kappa(y,(\tau+1)/2;\tau)$.
One can easily check that this answer agrees with the corresponding
Fukaya product. Note also that in fact this Massey product can
still be expressed via theta functions due to the identity
$$\kappa(y,\frac{\tau+1}{2})=
\frac{\theta'(\frac{\tau+1}{2})}{2\pi i \theta(y-\frac{\tau+1}{2})}$$

If we replace $\xi$ by another point on $E$ the Massey product will no
longer be univalued. However, if one chooses a splitting of the embedding
$\LL\ra F_y$ over some open subset $U\subset E$ and a trivialization of
$\LL|_U$ then replacing $r$ by this splitting we get a univalued
operation. The function $\kappa(y,x;\tau)$ appears as such operation
corresponding to the choice of a splitting over $E\setminus S$
coming from the trivialization of the pull-back of $F_y$ to $\C^*$.

Another example of well-defined Massey products that are expressed
in terms of the function $\kappa$ is the following.
Consider an extension
$$0\ra\LL\stackrel{a}{\ra} F\stackrel{b}{\ra} M\ra 0$$
where $M$ and $\LL$ are line bundles, $\deg M>0$, $\deg\LL=0$,
$\LL\not\simeq\O$.
Then as we have seen before for any section $s:\O\ra M$ the
triple product $m_3(s,c,a)$, where $c:M\ra\LL[1]$ is represented
by the above extension, is a lifting of $s$ to a section of $F$.
The corresponding Massey product is well-defined so we have
$MP(s,c,a)=m_3(s,c,a)$. The latter Fukaya product is of trapezoid
type since $\deg\LL=\deg\O=0$, so it can be expressed via $\kappa$. 

\subsection{Associativity constraint} Let us consider an example
of associativity constraint for Fukaya's $A_{\infty}$-category
of a torus involving triple products computed above.
 
Let us consider the following five lagrangians
in $\R^2/\Z^2$: $\La_1=(t,-t)$ with trivial connection,
$\La_2=(t,\a_2)$ with $2\pi i\b_2dx$, $\La_3=(-\a_1-\a_2,t)$ with
$2\pi i(\b_1+\b_2)dy$, $\La_4=(t,\a_1+\a_2+\a_3)$ with
$2\pi i(\b_1+\b_2+\b_3)dx$, $\La_5=(0,t)$ with trivial connection.
Here $\a_i,\b_i$ are real numbers, $\a_1+\a_2\not\in\Z$,
$\a_1+\a_3\not\in\Z$.
We choose liftings of $\La_i$ to objects in Fukaya
category in such a way that there is a non-zero morphism $a_i$ of
degree 0 from $\La_i$ to $\La_{i+1}$. We want to write the above
$A_{\infty}$-identity for these morphisms.
Note that all the $\Hom$-spaces between our objects are either
zero or one-dimensional, so the relevant compositions $m_2$ and $m_3$
are just numbers. Taking into account the fact that
$\La_2\cap\La_4=\La_3\cap\La_5=\emptyset$ the identity boils down to
\begin{equation}\label{axiom}
\begin{array}{l}
m_2(\La_1,\La_2,\La_3)m_3(\La_1,\La_3,\La_4,\La_5)+
m_2(\La_1,\La_4,\La_5)m_3(\La_1,\La_2,\La_3,\La_4)-\\
m_2(\La_1,\La_2,\La_5)m_3(\La_2,\La_3,\La_4,\La_5)=0
\end{array}
\end{equation}
where for example $m_2(\La_1,\La_2,\La_3)$ is the unique non-zero
$m_2$-composition of morphisms between $\La_1,\La_2,\La_3$, etc.
   
Now we can express all ingredients of (\ref{axiom}) in terms of
theta-functions, and functions $f$ and $g$ introduced above.
Namely, denoting $z_i=\a_i\tau+\b_i$ for $i=1,2,3$ we obtain
\begin{align*}
&m_2(\La_1,\La_2,\La_3)=\e(\a_1^2\frac{\tau}{2}+\a_1\b_1)\theta(z_1,\tau),\\
&m_2(\La_1,\La_4,\La_5)=\\
&\e((\a_1+\a_2+\a_3)^2\frac{\tau}{2}+
(\a_1+\a_2+\a_3)(\b_1+\b_2+\b_3))\theta(z_1+z_2+z_3,\tau),\\
&m_2(\La_1,\La_2,\La_5)=\e(\a_2^2\frac{\tau}{2}+\a_2\b_2)
\theta(z_2,\tau),\\
&m_3(\La_1,\La_3,\La_4,\La_5)=\\
&\e((\frac{\a_1+\a_2}{2}+\a_3)(\a_1+\a_2)\tau+
(\a_1+\a_2)(\b_1+\b_2+\b_3)+\a_3(\b_1+\b_2))g(z_3,z_1+z_2;\tau),\\
&m_3(\La_1,\La_2,\La_3,\La_4)=\e((\frac{\a_1^2}{2}-\frac{\a_3^2}{2})\tau+
\a_1\b_1-\a_3\b_3)g(-z_3,z_1+z_3;\tau),\\
&m_3(\La_2,\La_3,\La_4,\La_5)=\\
&\e((\a_1+\a_2)(\a_1+\a_3)\tau+
(\a_1+\a_2)(\b_1+\b_3)+(\a_1+\a_3)(\b_1+\b_2))f(z_1+z_2,z_1+z_3;\tau).
\end{align*}
Substituting this into (\ref{axiom}) and deleting similar terms we
obtain the following identity
\begin{equation}\label{fg}
\begin{array}{l}
\theta(z_1,\tau)g(z_3,z_1+z_2;\tau)+
\theta(z_1+z_2+z_3,\tau)g(-z_3,z_1+z_3;\tau)=\\
\theta(z_2,\tau)f(z_1+z_2,z_1+z_3;\tau).
\end{array}
\end{equation}
By (\ref{g0}) this implies the following identity between meromorphic
functions of $\C^3$:
\begin{equation}\label{identity1}
\e(y)\theta(y+z,\tau)\kappa(y,z-x;\tau)-
\e(-x)\theta(x-z,\tau)\kappa(-x,y+z;\tau)=
\theta(z,\tau)f(x,y;\tau)
\end{equation}
where we have put $x=z_1+z_2$, $y=z_1+z_3$, $z=z_2$.

We used the fact that Fukaya composition satisfies axioms
of $A_{\infty}$-category to derive the above identity.
However, it can be also proved in a straightforward way
comparing residues of both sides at poles and
using the difference equation for $\kappa$.
It appears in  a  slightly  different  form in
Halphen's  book  \cite{H}  (p.481, formula (45) and the next
one).

\subsection{More Fukaya products}
Let us consider another example of triple Fukaya products where
none of the four lines are parallel. Namely, let
$\La_1=(\a_1+t,-\a_1+t)$, $\La_2=(-\a_2,t)$, $\La_3=(t,0)$,
$\La_4=(t,-t)$. Then we can choose (essentially uniquely) lifts of
the corresponding lagrangain circles in $\R^2/\Z^2$ to objects
of Fukaya category in such a way that there is a non-zero morphism
of degree zero from $\La_i$ to $\La_{i+1}$. All the Hom-spaces
$\Hom(\La_i,\La_{i+1})$ are 1-dimensional, however, $\Hom^{-1}(\La_1,\La_4)$
is 2-dimensional since the corresponding circles have 2 points of
intersection. Let us consider the component of the triple product
corresponding to the point $(\a_1,-\a_1)\in\La_1\cap\La_4$.
Then as before its value is expressed via certain holomorphic function
of $z_1=\a_1\tau$ and $z_2=\a_2\tau$ (one could add some monodromies
in the above picture to get all values of complex
variables $z_1$ and $z_2$). This function has the following form
\begin{align*}
&h(z_1,z_2;\tau)=\\
&\sum_{(m+\a(z_1))(n+\a(z_2))>0}\sign(m+\a(z_1))
\e\left(\frac{\tau}{2}\left(2m^2+4mn+n^2\right)+2(m+n)z_1+(2m+n)z_2\right).
\end{align*}
It is holomorphic in the region $\a(z_i)\not\in\Z$ and can be considered
as a section of a line bundle on $E^2$ over the corresponding open subset of
$E^2$ because of the following quasi-periodicity:
\begin{eqnarray}\label{hqp1}
h(z_1+m+\tau,z_2;\tau)=\e(-\tau-2z_1-2z_2)h(z_1,z_2;\tau),\\
\label{hqp2}
h(z_1,z_2+m+\tau;\tau)=\e(-\tau/2-2z_1-z_2)h(z_1,z_2;\tau).
\end{eqnarray}
Let us denote by $h_0$ the restriction of $h$ to the region
$0<\a(z_i)<1$ for $i=1,2$. Then it is easy to see that $h_0$ extends to
a holomorphic function on $\C^2$ (the series converges absolutely).

As before we can translate certain $A_{\infty}$-associativity axiom into an
identity involving $h(z_1,z_2)$. More precisely, we can consider five lines:
two with slope 0, one with slope 1, one with slope -1, and one with slope
$\infty$. Their relative position is described by three parameters. Adding
monodromies we get the following identity with three complex variables
$z_1,z_2,z_3$:
\begin{equation}
\begin{array}{l}
\theta(2z_1+z_3,\tau)h(z_1+z_3,-z_1+z_2-z_3;\tau)+\theta(z_1+z_2+z_3,\tau)
h(-z_1-z_3,z_3;\tau)=\\
=\theta(2z_2,2\tau)g(-z_1+z_2-z_3,-z_1-z_2;\tau)+
\theta(2z_1,2\tau)g(z_3,z_1+z_2)
\end{array}
\end{equation}
Using difference equations for $h$ and $g$
we can further transform this into the following
identity involving $h_0$ and $\kappa$:
\begin{equation}\label{identity2}
\begin{array}{l}
\th(2x+y,\tau)h_0(x,z;\tau)-\th(2x+z,\tau)h_0(x,y;\tau)=\\
=\th(2(x+z),2\tau)\kappa(-2x-y-z,2x+y+\tau;\tau)-\\
-\th(2(x+y),2\tau)\kappa(-2x-y-z,2x+z+\tau;\tau)
\end{array}
\end{equation}

Substituting $y=-x$ in this identity we can express $h_0(x,z;\tau)$ via
$h_0(x,-x;\tau)$ and functions $\kappa$ and $\theta$. We claim that
$h_0(x):=h_0(x,-x;\tau)$ can also be expressed as a rational function of
$\kappa$ and
$\theta$. Indeed, we have
$$h_0(x)=\sum_{(m+1/2)(n+1/2)>0}\sign(m+1/2)
\e\left(\frac{\tau}{2}\left(2m^2+4mn+n^2\right)+nx\right).$$
Hence,
\begin{align*}
&h_0(x+\tau)=\e(\frac{\tau}{2}+x)\times\\
&\sum_{(m+1/2)(n+1/2)>0}\sign(m+1/2)
\e\left(\frac{\tau}{2}\left(2(m+1)^2+4(m+1)(n-1)+(n-1)^2\right)+(n-1)x\right)
=\\
&=\e(\frac{\tau}{2}+x)\left(h_0(x)-\th(x,\tau)\right)+\th(0,2\tau).
\end{align*}
Let us denote $\psi(x)=\th(x-(\tau+1)/2,\tau)h_0(x)$. Then $\psi$ satisfies the
equation
$$\psi(x+\tau)=\e(\frac{\tau+1}{2})\psi(x)+\e(\frac{\tau}{2})
\th(x,\tau)\th(x-\frac{\tau+1}{2},\tau)+
\th(0,2\tau)\th(x+\frac{\tau+1}{2},\tau).
$$
It is easy to see that this equation has the unique holomorphic solution,
namely,
\begin{align*}
&\psi(x)=\th(0,2\tau)\kappa(\frac{\tau+1}{2},x+\frac{\tau+1}{2};\tau)+\\
&\th(\frac{\tau+1}{2},2\tau)\cdot\left(\e(\frac{\tau}{2})
\kappa(\frac{\tau+1}{2},2x-\frac{\tau+1}{2};2\tau)-
\e(x-\frac{\tau}{2})\kappa(-\frac{\tau+1}{2},2x+\frac{\tau+1}{2};2\tau)\right)
\end{align*}

\subsection{Generic case}
Let us consider four lines
$L_i(y_i)=\{(\frac{y_i+t}{\la_i},t)\}$, $i=1,2,3,4$, where
the slopes $\la_i$ are distinct rational numbers, $y_i$
are fixed real numbers.
Let us denote by $\ov{L}_i(y_i)$ the image of $L_i$ in $\R^2/\Z^2$.
The intersection of $\ov{L}_i(y_i)$ and $\ov{L}_j(y_j)$
consists of the following set of points:
$$e_{a,b}(y_i,y_j)=(y_{ij}, y'_{ij})+\frac{a\la_j+b}{\la_j-\la_i}(1,\la_i)
\mod\Z^2$$
where $a,b$ are integers, $y_{ij}=\frac{y_j-y_i}{\la_j-\la_i}$,
$y'_{ij}=\frac{\la_iy_j-\la_jy_i}{\la_j-\la_i}$.
Note that if $a$ and $b$ are integers then
$\ov{L}_i(y_i+a\la_i+b)=\ov{L}(y_i)$. Furthermore, we have
$$e_{a,b}(y_i,y_j)=e_{0,0}(y_i-a\la_i-b,y_j)=
e_{0,0}(y_i,y_j+a\la_j+b).$$

For every $i<j$ and an intersection point
$p\in\ov{L}_i(y_i)\cap\ov{L}_j(y_j)$
we   denote   by   $[p]$  the  corresponding  morphism  from
$\ov{L}_i(y_i)$ to $\ov{L}_j(y_j)$ in the  Fukaya
category. Also let us denote $e_{ij}=e_{0,0}(y_i,y_j)$.
Note that $\deg[e_{a,b}(y_i,y_j)]=\deg(i,j)$ where
$\deg(i,j)=0$ if $\la_i<\la_j$ and
$\deg(i,j)=1$ otherwise.
We want to compute the Fukaya product
$m_3([e_{12}],[e_{23}],[e_{34}])$, so we assume that
\begin{equation}\label{degree}
\sum_{i=1}^3\deg(i,i+1)=\deg(1,4)+1
\end{equation}
(otherwise this product is zero).

For every rational number $\la$ let us denote by $I_{\la}\sub\Z$
the subset of $n\in\Z$ such that $n\la\in\Z$.
To each quadruple of distinct rational numbers $\la_1,\ldots,\la_4$
we associate a lattice $\La=\La(\la_1,\la_2,\la_3,\la_4)$ and a sublattice
$\La^+=\La^+(\la_1,\la_2,\la_3,\la_4)$ as follows:
\begin{align*}
&\La(\la_1,\la_2,\la_3,\la_4)=\\
&=\{\ov{n}=(n_1,n_2,n_3,n_4)\in\Q^4\ |\
\sum_{i=1}^4 n_i=0, \sum_{i=1}^4 \la_in_i=0;
n_2\in I_{\la_2}, n_3\in I_{\la_3}\},\\
&\La^+(\la_1,\la_2,\la_3,\la_4)=\{\ov{n}\in
\La(\la_1,\la_2,\la_3,\la_4)\ |\ n_1\in I_{\la_1}\}.
\end{align*}
Equivalently, the sublattice $\La^+$ is distinguished in $\La$ by the
condition $n_4\in I_{\la_4}$.
Consider the following quadratic form on $\La\ot\R$:
$$Q(x)=Q_{\la_1,\la_2,\la_3,\la_4}(x)=(\la_3-\la_4)x_3x_4+(\la_1-\la_2)x_1x_2.
$$
Clearly, $Q$ takes integer values on $\La^+$.
However, it is in general indefinite.
Let  $C\sub  \La\ot\R\simeq\R^2$ be the subset
consisting of $x$ such that
$(\la_{i-1}-\la_i)x_ix_{i-1}>0$ for all $i=1,2,3,4$
(where   we   set   $x_0=x_4$).
Then clearly $Q(x)>0$ for any $x\in C$.
Note  that  our  assumption
(\ref{degree}) implies that $C$ is non-empty.
It is easy to see that
two of the four inequalities defining $C$ are redundant, so
$C$  is  always a region in the plane bounded by 2 lines.
Let  $C=C^+\cup C^-$ be a
decomposition of $C$ into connected components,
$\eps:C\ra  \pm 1$ be the function assigning $1$ (resp.
$-1$) to points of $C^+$ (resp. $C^-$).
Note    that    for    any   quadruple   of   real   numbers
$y=(y_1,y_2,y_3,y_4)$ the vector
$$v(y)=(y_{14}-y_{12},y_{12}-y_{23},y_{23}-y_{34},
y_{34}-y_{14})$$
(where $y_{ij}=\frac{y_j-y_i}{\la_j-\la_i}$)
belongs to the subspace $\La\ot\R\sub\R^4$.
Now for $\ov{n}_0\in\La\ot\Q$ and $z\in\C^4$ we set
\begin{align*}
&F_{\la_1,\la_2,\la_3,\la_4;\ov{n}_0}(z)=\\
&=\sum_{\ov{n}\in(\La^+ +\ov{n}_0)\cap(C-v(\a(z)))}
\eps(\ov{n}+v(\a(z)))\cdot
\e(\frac{\tau}{2} Q(\ov{n})+n_1z_1+n_2z_2+n_3z_3+n_4z_4))
\end{align*}
where  $\tau\cdot\a$  is  the  first  projection  of  the  direct sum
decomposition $\C^4= \tau\R^4\oplus \R^4$. For
$z$ varying in some open subset of $\C^4$ (in fact, in
the complement to codimension one analytic subset), this
is a holomorphic function of $z$.
Note that linear relations between $n_i$ imply that
$$n_1z_1+n_2z_2+n_3z_3+n_4z_4=(\la_1-\la_2)(z_{14}-z_{12})n_2+
(\la_1-\la_3)(z_{14}-z_{13})n_3,$$
where $z_{ij}=\frac{z_j-z_i}{\la_j-\la_i}$.
On the other hand, the vector $v(\a(z))$ is determined by its
first two components $\a(z_{14}-z_{12})$ and
$\a(z_{12}-z_{23})=\frac{\la_3-\la_1}{\la_3-\la_2}\a(z_{12}-z_{13})$.
Thus, $F$ actually depends on two holomorphic variables
$z_{13}-z_{12}$ and $z_{14}-z_{12}$.

\begin{prop}                   Assume                   that
$\sum_{i=1}^3\deg[e_{i,i+1}]=\deg[e_{1,4}]+1$. Then one has
\begin{eqnarray}\label{m3}
m_3([e_{12}],[e_{23}],[e_{34}])=\pm
\e(\frac{\tau}{2}\Delta(y_1,y_2,y_3,y_4))\cdot\\
\sum_{\ov{k}\in\La/\La^+}
F_{\la_1,\la_2,\la_3,\la_4;\ov{k}}(\tau y)
[e_{-k_2-k_3,\la_2k_2+\la_3k_3}(y_1,y_4)].\nonumber
\end{eqnarray}
where
$$\De(y_1,y_2,y_3,y_4)=\det\left(\matrix   y_{34}-y_{12}        &
y_{23}-y_{14}           \\          y'_{34}-y'_{12}          &
y'_{23}-y'_{14}\endmatrix\right)$$
\end{prop}

\Pf. The idea is to parametrize the quadrangles contributing to $m_3$
by elements of $\La$. Namely, let $p_{i, i+1}$, $i=1,\ldots, 4$
be the vertices of such a quadrangle (where the edge between
$p_{i-1,i}$ and $p_{i,i+1}$ belongs to $\ov{L}_i(y_i)$ modulo $\Z^2$).
Then denoting the difference between the first coordinates of
$p_{i-1,i}$ and $p_{i,i+1}$ by $n_i$ we find that $(n_1,n_2,n_3,n_4)$
belongs to $\La(\la_1,\la_2,\la_3,\la_4)$. The condition
$n+v(y)\in C$ is equivalent to the requirement that the lagrangians come
in the correct order when one goes clockwise along the quadrangle.
\ed

In  the  above  proposition  we   consider   $\ov{L}_i(y_i)$
equipped  with  trivial  local systems. If we add non-trivial
connections   along   $\ov{L}_i(y_i)$   the formula (\ref{m3}) will change
by adding some linear combinations of these connections to the
arguments of the function $F$.
           
For five lines $L(y_i)$, $i=1,\ldots,5$ we can consider the
identity obtained by comparing the coefficients
with $[e_{15}]$ in the $A_{\infty}$-constraint (\ref{a-infty})
for $a_1=[e_{12}]$, $a_2=[e_{23}]$, $a_3=[e_{34}]$ and
$a_4=[e_{45}]$.

Before writing the corresponding identity let us introduce
some more notation.
Let  $\la_1,\la_2,\la_3$  be  a triple of distinct rational
numbers. We denote
$$I_{\la_1,\la_2,\la_3}=I_{\la_2}\cap
\frac{\la_3-\la_1}{\la_3-\la_2}I_{\la_1}.$$
One              easily              checks             that
$I_{\la_1,\la_2,\la_3}=I_{\la_3,\la_2,\la_1}$.
Assume that
\begin{equation}\label{degree2}
\deg(1,2)+\deg(2,3)=\deg(1,3).
\end{equation}
Then similarly to the previous proposition one checks that
\begin{align*}
&m_2([e_{12}],[e_{23}])=\\
&\e(-\frac{\tau}{2}\De(y_1,y_2,y_3))\cdot\sum_{n_0\in
I_{\la_2}/I_{\la_1,\la_2,\la_3}}\th_{\la_1,\la_2,\la_3;n_0}(
\tau(y_1,y_2,y_3))[e_{n_0,-\la_2n_0}(y_1,y_3)]
\end{align*}
where
$$\De(y_1,y_2,y_3)=\det\left(\matrix       y_{23}-y_{12}        &
y_{13}-y_{12}           \\          y'_{23}-y'_{12}          &
y'_{13}-y'_{12}\endmatrix\right),$$
\begin{align*}
&\th_{\la_1,\la_2,\la_3;n_0}(z_1,z_2,z_3)=\\
&=\sum_{n\in
I_{\la_1,\la_2,\la_3}}\e\left(
\frac{(\la_3-\la_2)(\la_2-\la_1)\tau}{2(\la_3-\la_1)}(n+n_0)^2+
(n+n_0)(z_{23}-z_{12})\right).
\end{align*}
Note that  the condition (\ref{degree2}) implies that
$(\la_3-\la_2)(\la_2-\la_1)(\la_3-\la_1)^{-1}>0$ so the above
series converges.

Now  the  $A_{\infty}$-constraint  leads  to  the  following
identity.

\begin{align*}
&\eps_1\cdot\sum_{\ov{n}\in\La_{2345}/\La^+_{2345},n_2\in
I_2+\frac{\la_5-\la_1}{\la_5-\la_2}I_1}
F_{2345,\ov{n}}(z)\th_{125,n''_2}(z)+\\
&\eps_2\cdot\sum_{\ov{n}\in\La_{1234}/\La^+_{1234},n_4\in
I_4+\frac{\la_5-\la_1}{\la_4-\la_1}I_5}
F_{1234,\ov{n}}(z)\th_{145,n''_4}(z)+\\
&\eps_3\cdot\sum_{k\in I_2/I_{123},
k\in \frac{\la_1-\la_3}{\la_1-\la_2}I_3
+\frac{\la_1-\la_4}{\la_1-\la_2}I_4+\frac{\la_1-\la_5}{\la_1-
\la_2}I_5}
F_{1345,k''u''+k'''u'''}(z)\th_{123,k}(z)+\\
&\eps_4\cdot\sum_{k\in I_4/I_{345},
k\in \frac{\la_5-\la_1}{\la_5-\la_4}I_1
+\frac{\la_5-\la_2}{\la_5-\la_4}I_2+\frac{\la_5-\la_3}{\la_5-
\la_4}I_3}
F_{1235,k'v'+k''v''}(z)\th_{345,k}(z)+\\
&\eps_5\cdot\sum_{k\in I_3/I_{234},
k\in \frac{\la_5-\la_1}{\la_5-\la_3}I_1
+\frac{\la_5-\la_2}{\la_5-\la_3}I_2+\frac{\la_5-\la_4}{\la_5-
\la_3}I_4}
F_{1245,k'w'+k''w''+k'''w'''}(z)\th_{234,k}(z)=0.
\end{align*}

Here $z$ is in $\C^5$ minus some analytic subset of codimension 1,
however, all the functions in this identity actually depend on
three holomorphic variables: $z_{13}-z_{12}$, $z_{14}-z_{12}$ and
$z_{15}-z_{12}$ (where as before $z_{ij}=(z_j-z_i)/(\la_j-\la_i)$).
For every $1\le i<j<k<l\le 5$ we denote
$F_{ijkl,\ov{n}}(z)=F_{\la_i,\la_j,\la_k,\la_l;\ov{n}}(
z_i,z_j,z_k,z_l)$. Similarly,
$\th_{ijk,n}(z)=\th_{\la_i,\la_j,\la_k;n}(z_i,z_j,z_k)$.
Also we write for simplicity $I_i=I_{\la_i}$,
$\La_{ijkl}=\La(\la_i,\la_j,\la_k,\la_l)$, etc.
The elements of $\La_{ijkl}\ot\Q$ are denoted by
$\ov{n}=(n_i,n_j,n_k,n_l)$.
The multiple $\eps_i$ before each term is $\pm 1$ if
the conditions on the degrees are satisfied (we will specify the sign later),
and $0$ otherwise.
For example, $\eps_1=0$ unless
$\deg(1,2)+\deg(2,5)=\deg(1,5)$ and
$\deg(2,3)+\deg(3,4)+\deg(4,5)=\deg(2,5)+1$.
In the first two terms of the identity we denote
$n_2=n_2'+n_2''$ (resp. $n_4=n_4'+n_4''$) with
respect to the inclusion
$n_2\in I_2+\frac{\la_5-\la_1}{\la_5-\la_2}I_1$
(resp. $n_4\in I_4+\frac{\la_5-\la_1}{\la_4-\la_1}I_5$).
Similarly, in the last three terms the decomposition
$k=k'+k''+k'''$ corresponds to the inclusion of $k$
into sum of three  ideals.  Note  that  although  $k$  is  a
representative of a coset,
the condition that $k$ belongs to the sum of three ideals
is well-defined. In 3rd and 4th term this is clear, while in
the 5th term the inclusion
$$I_{234}\sub\frac{\la_5-\la_2}{\la_5-\la_3}I_2+
\frac{\la_5-\la_4}{\la_5-\la_3}I_4$$
follows from the identity
$$k=\frac{\la_5-\la_2}{\la_5-\la_3}\cdot(\frac{\la_4-\la_3}
{\la_4-\la_2}k)+\frac{\la_5-\la_4}{\la_5-\la_3}\cdot
(\frac{\la_3-\la_2}{\la_4-\la_2}k).$$
Finally, we denoted
\begin{align*}
&u'=(\frac{(\la_4-\la_3)(\la_2-\la_1)}{(\la_4-\la_1)(\la_3-\la_1)},
\frac{\la_1-\la_2}{\la_3-\la_1},\frac{\la_1-\la_2}{\la_1-\la_4},0),\\
&u''=(\frac{(\la_5-\la_3)(\la_2-\la_1)}{(\la_5-\la_1)(\la_3-\la_1)},
\frac{\la_1-\la_2}{\la_3-\la_1},0,\frac{\la_2-\la_1}{\la_5-\la_1}),\\
&v'=(\frac{\la_5-\la_4}{\la_5-\la_1},0,\frac{\la_5-\la_4}{\la_3-\la_5},
\frac{(\la_5-\la_4)(\la_1-\la_3)}{(\la_5-\la_1)(\la_3-\la_5)}),\\
&v''=(0,\frac{\la_5-\la_4}{\la_5-\la_2},\frac{\la_5-\la_4}{\la_3-\la_5},
\frac{(\la_5-\la_4)(\la_2-\la_3)}{(\la_5-\la_2)(\la_3-\la_5)}),\\
&w'=(\frac{\la_5-\la_3}{\la_5-\la_1},\frac{\la_4-\la_3}{\la_2-\la_4},
\frac{\la_3-\la_2}{\la_2-\la_4},\frac{\la_3-\la_1}{\la_5-\la_1}),\\
&w''=(0,\frac{(\la_5-\la_4)(\la_2-\la_3)}{(\la_5-\la_2)(\la_2-\la_4)},
\frac{\la_3-\la_2}{\la_2-\la_4},\frac{\la_3-\la_2}{\la_5-\la_2}),\\
&w'''=(0,\frac{\la_4-\la_3}{\la_2-\la_4},
\frac{(\la_5-\la_2)(\la_3-\la_4)}{(\la_2-\la_4)(\la_5-\la_4)},
\frac{\la_3-\la_4}{\la_5-\la_4})
\end{align*}

Note that there is an ambiguity of sign in the definition of
$F_{ijkl}$. The claim is that for
every given configuration of $\deg(e_i,e_j)$ there exists
a choice of signs in $F_{ijkl}$ and $\eps_i$ which makes the
above identity true. For example, let us assume that
$\la_3<\la_1<\la_4<\la_2<\la_5$. Then all $\eps_i$ are non-zero.
Let us choose the positive components $C^+$ of the cones $C_{ijkl}$
entering in the definition of $F_{ijkl}$ as follows:
\begin{align*}
&C^+_{2345}: n_2>0, n_3>0, n_4<0, n_5>0,\\
&C^+_{1234}: n_1>0, n_2<0, n_3<0, n_4>0,\\
&C^+_{1345}: n_1>0, n_3>0, n_4<0, n_5>0,\\
&C^+_{1235}: n_1>0, n_2<0, n_3<0, n_5>0,\\
&C^+_{1245}: n_1>0, n_2<0, n_4<0, n_5>0.
\end{align*}
Then the signs should be chosen as follows: $\eps_1=\eps_2=\eps_5=1$,
$\eps_3=\eps_4=-1$. To see this we note that for purely imaginary $\tau$
we can represent each of five terms of our identity in the form
$$\sum_{\ov{n},m} \eps(\ov{n}) c_{\ov{n},m}\e(\sum_{i=1}^5 l_i(\ov{n},m)z_i)$$
where the sum is taken over a coset for a lattice in
$(\La_{ijkl}\otimes\Q)\times\Q$, $l_i$ are some rational linear functions
of $(\ov{n},m)$, and the coefficients $c_{\ov{n},m}$ are positive.
It is easy to see that $l_i$ are linearly independent, so there are
no cancellations in the above Fourier series. Now we consider
a term corresponding to $\ov{n}\in C^+_{2345}$ and $m>>0$ (resp. $m<<0$)
and find out
that the only term it can cancel out with has $\ov{n'}\in C^+_{1235}$
(resp. $\ov{n'}\in C^-_{1245}$), hence, $\eps_4=-\eps_1$ (resp.
$\eps_5=\eps_1$). On the other hand, a term corresponding to
$\ov{n}\in C^+_{1234}$ and $m>>0$ (resp. $m<<0$) can cancel out only with
the term which has $\ov{n'}\in C^+_{1345}$ (resp. $\ov{n'}\in C^-_{1245}$),
hence, $\eps_3=-\eps_2$ (resp. $\eps_5=\eps_2$).


\begin{thebibliography}{99}
\bibitem{Ap} M.~P.~Appell, {\it Sur le fonctions doublement
periodique de troisieme espece},  Annales  scientifiques  de
l'\'Ecole Normale Sup\'erieure, 3e s\'erie, t.I, p.135,
t.II, p.9, t.III, p.9 (1884--1886).
\bibitem{GR} G.~Gasper, M.~Rahman, {\it Basic hypergeometric series},
Cambridge University Press, 1990.
\bibitem{GM}  S.~Gelfand, Yu.~Manin, {\it Methods of
homological algebra}. Springer-Verlag, 1996.
\bibitem{F} K.~Fukaya, {\it Morse Homotopy, $A^\infty$-Category, and
Floer Homologies}, in {\it The Proceedings of the
1993 GARC Workshop on Geometry and Topology,} H.~J.~Kim, ed.,
Seoul National University.
\bibitem{F-rec} K.~Fukaya, {\it Mirror symmetry  of  Abelian
variety and Multi Theta functions}, preprint, 1998.
\bibitem{GZ} L. G\"ottsche, D. Zagier, {\it Jacobi forms and the structure
of Donaldson invariants for 4-manifolds with $b_+=1$}. Selecta Mathematica,
4 (1998), 69--115.
\bibitem{H}  G.-H.~Halphen,  {\it  Trait\'e  des  fonctions
elliptiques}, I. Paris, 1886.
\bibitem{Kon} M. Kontsevich, {\it Homological algebra of mirror symmetry},
Proceedings of ICM (Z\"urich, 1994), 120--139. Birkh\"auser, Basel, 1995.
\bibitem{Kr} L. Kronecker, {\it Zur theorie der elliptischen
functionen} (1881), in {\it Leopold Kronecker's Werke}, vol. IV, 311--318.
Chelsea Pub. Co., 1968.
\bibitem{PZ} A. Polishchuk, E. Zaslow, {\it Categorical mirror
symmetry: the elliptic curve}, Adv. in Theor. and Math. Physics 
2 (1998), 443--470.
\bibitem{Weil} A. Weil, {\it Elliptic functions according to
Eisenstein and Kronecker}. Springer-Verlag, 1976.
\bibitem{Z} D.~Zagier, {\it Periods of modular forms and
Jacobi theta functions}, Invent. Math. 104 (1991), 449--465.
\end{thebibliography}
\end{document}